\documentclass[twoside]{article}
\usepackage{graphicx}
\usepackage{amsmath, amsthm, amssymb, mathrsfs}

\newtheorem{theorem}{Theorem}
\newtheorem{remark}{Remark}
\newtheorem{lemma}{Lemma}
\newtheorem{corollary}{Corollary}
\newfont{\mi}{cmti9}
\newfont{\m}{cmr8}
\newfont{\ms}{cmsl8}
\newfont{\autor}{cmcsc10}
\begin{document}
\rule[3mm]{128mm}{0mm}
\vspace*{-16mm}

{\footnotesize K\,Y\,B\,E\,R\,N\,E\,T\,I\,K\,A\, ---
\,V\,O\,L\,U\,M\,E\, {\it 4\,5\,} (\,2\,0\,0\,9\,)\,,\,
N\,U\,M\,B\,E\,R\, X\,,\,\  P\,A\,G\,E\,S\, \,x\,x\,x -- x\,x\,x}\\
\rule[3mm]{128mm}{0.2mm}

\vspace*{11mm}

{\large\bf \noindent OPTIMAL SEQUENTIAL PROCEDURES WITH BAYES
DECISION RULES}

\vspace*{8mm}

{\autor \indent Andrey Novikov}\vspace*{23mm}\\
\small In this article, a general problem of sequential statistical
inference for general discrete-time stochastic processes is
considered. The problem is to minimize an average sample number
given that Bayesian risk due to incorrect decision does not exceed
some given bound. We characterize the form of optimal sequential
stopping rules in this problem. In particular, we have a
characterization of the form of optimal sequential decision
procedures when the Bayesian risk includes both the loss due to
incorrect decision and the cost of observations.

\smallskip\par
\noindent {\sl Keywords:}\,
\begin{minipage}[t]{112mm}
 sequential analysis, discrete-time stochastic
process, dependent observations, statistical decision problem, Bayes
decision, randomized stopping time, optimal stopping rule,
existence and uniqueness of optimal sequential decision procedure
\end{minipage} \smallskip

\normalsize

\section{\normalsize INTRODUCTION}\label{s1}

Let $X_1,X_2,\dots, X_n, \dots$ be a discrete-time stochastic
process, whose distribution depends on   an unknown parameter
$\theta$, $\theta\in\Theta$. In this article, we consider a general problem of
sequential statistical decision making based on the observations of
this process.

Let us  suppose that for any $n=1,2,\dots,$ the vector $(X_1, X_2,
\dots , X_n)$ has a probability ``density'' function
\begin{equation}\label{0}
 f_\theta^n=f_\theta^{n}(x_1, x_2, \dots, x_n)
 \end{equation}
 (Radon-Nikodym derivative of its distribution) with respect to
 a
 product-measure $$\mu^n=\mu\otimes \mu\otimes \dots\otimes\mu,$$
 with some $\sigma$-finite measure $\mu$ on the respective space.
 As usual in the Bayesian context, we suppose that $f_\theta^{n}(x_1, x_2, \dots,
 x_n)$ is measurable with respect to $(\theta,x_1,\dots,x_n)$, for
 any $n=1,2,\dots$.

Let us  define a {\em sequential statistical procedure} as a pair
$(\psi,\delta)$, being $\psi$ a (randomized) {\em stopping rule},
$$\psi=\left(\psi_1,\psi_2, \dots ,\psi_{n},\dots\right),$$
and $\delta$ a {\em decision rule}
$$
\delta=\left(\delta_1,\delta_2, \dots ,\delta_{n},\dots\right),
$$
supposing that
$$\psi_n=\psi_n(x_1, x_2,\dots, x_n)$$  and
$$
\delta_n=\delta_n(x_1, x_2,\dots, x_n)
$$
are measurable functions,  $\psi_n(x_1,\dots,x_n)\in [0, 1]$, $\delta_n(x_1,\dots,x_n)\in\mathscr D$ (a decision space), for
any observation vector $(x_1,\dots, x_n)$,  for any $n=1,2,\dots$
(see, for example, \cite{Wald}, \cite{DeGroot}, \cite{Ferguson}, \cite{Berger},
\cite{Ghosh}).

The interpretation of  these elements is as follows.

The value of  $\psi_n(x_1,\dots,x_n)$ is interpreted  as the
conditional probability {\em to stop and proceed to decision
making}, given that that we came to stage $n$ of the experiment and
that the observations up to stage $n$ were $(x_1, x_2, \dots, x_n).$
If there is no stop, the experiments continues to the next stage and
an additional observation $x_{n+1}$ is taken. Then the rule
$\psi_{n+1}$ is applied to $x_1,\dots, x_n,x_{n+1}$ in the same way
as as above, etc., until the experiment eventually stops.

When the experiments stops at stage $n$, being $(x_1,\dots,x_n)$ the
data vector observed, the decision specified by $\delta_n(x_1,\dots,x_n)$ is taken, and the sequential 
statistical experiment stops.

The stopping rule $\psi$ generates, by the above process, a random
variable $\tau_\psi$ (randomized {\em stopping time}), which may be  defined as follows.
Let $U_1,U_2\dots,U_n,\dots$ be a sequence of independent and identically distributed (i.i.d.)
random variables  uniformly distributed on $[0,1]$ (randomization variables), such that the
process $(U_1,U_2,\dots)$ is independent of the process of observations $(X_1,X_2,\dots)$.
Then let us say that $\tau_\psi=n$ if, and only if,
$$U_1>\psi_1(X_1),\; \dots, \;U_{n-1}>\psi_{n-1}(X_1,\dots,X_{n-1}),\;
 \mbox{and}\; U_n\leq \psi_{n}(X_1,\dots X_{n}),$$ $n=1,2,\dots$.

It is easy to see that the  distribution of $\tau_\psi$ is
given by
\begin{equation}\label{1}
P_\theta(\tau_\psi=n)=E_\theta (1-\psi_1)(1-\psi_2)\dots
(1-\psi_{n-1})\psi_n,\quad n=1,2,\dots.
\end{equation}

In (\ref{1}),  $\psi_n$ stands for
$\psi_n(X_1,\dots,X_n),$ unlike its previous definition as $\psi_n=\psi_n(x_1,\dots,x_n)$.
We use this ``duality" throughout the paper, applying, for any $F_n=F_n(x_1,\dots,x_n)$ or
$F_n=F_n(X_1,\dots X_n)$ the following general rule:
 when $F_n$ is under the probability or expectation
sign, it is $F_n(X_1,\dots, X_n)$, otherwise it is $F_n(x_1,\dots,
x_n)$.

Let $w(\theta, d)$ be a non-negative loss function (measurable with
 respect to $(\theta,d)$, $\theta\in\Theta$, $d\in\mathscr D$) and $\pi_1$
any probability measure on $\Theta$. We define  the average loss of the
sequential statistical procedure $(\psi,\delta)$ as
\begin{equation}\label{203}
W(\psi,\delta)=\sum_{n=1}^\infty \int
\left[E_\theta(1-\psi_1)\dots(1-\psi_{n-1})\psi_nw(\theta,\delta_n)\right]d\pi_1(\theta).
\end{equation}
and its {\em average sample number}, given $\theta$, as
\begin{equation}\label{8aa}
N(\theta;\psi)=E_\theta\tau_\psi
\end{equation}
(we suppose that $N(\theta;\psi)=\infty$ if $\sum_{n=1}^\infty
P_\theta(\tau_\psi=n)<1$ in (\ref{1})).

 Let us also define its "weighted" value
\begin{equation}\label{201} N(\psi)=\int
N(\theta;\psi)d\pi_2(\theta),
\end{equation}
where $\pi_2$ is some probability  measure on $\Theta$, giving ``weights'' to the
particular values of $\theta$.

 Our main
goal is minimizing $N(\psi)$ over all sequential decision procedures $(\psi,\delta)$ subject
to \begin{equation}\label{8aaa}W(\psi,\delta)\leq w,
\end{equation}
where $w$ is some positive constant, supposing that $\pi_1$ in
(\ref{203}) and $\pi_2$ in (\ref{201}) are, generally speaking, two
{\em different} probability measures. We only consider the cases when there exist procedures $(\psi,\delta)$ satisfying (\ref{8aaa}).

Sometimes it is necessary to put  the  risk under control in a more
detailed way. Let $\Theta_1, \dots, \Theta_k$ be some subsets of the
parametric space such that $\Theta_i\bigcap\Theta_j=\emptyset$ if
$i\not=j$, $i,j=1,\dots, k$. Then, instead of (\ref{8aaa}), we may
want to guarantee that
\begin{equation}\label{204}
W_i(\psi,\delta)=\sum_{n=1}^\infty\int_{\Theta_i}E_\theta(1-\psi_{i-1})\dots(1-\psi_{n-1})\psi_nw(\theta,\delta_n)d\pi_1(\theta)\leq
w_i,
\end{equation}
 with some $w_i>0$, for any $i=1,\dots,k$, when minimizing $N(\psi)$.

To advocate restricting the sequential procedures by (\ref{204}),
let us see a particular case of hypothesis testing.

Let $H_1:\theta=\theta_1$ and $H_2:\theta=\theta_2$ be two simple
hypotheses about the parameter value, and let
$$w(\theta,d)=\begin{cases}1& \mbox{if}\quad \theta=\theta_1\; \mbox{and}\; d=2,\cr
1& \mbox{if}\quad \theta=\theta_2\;\mbox\; \mbox{and}\; d=1,\cr 0&
\mbox{otherwise},\end{cases}$$ and $\pi_1(\{\theta_1\})=\pi$,
$\pi_1(\{\theta_2\})=1-\pi$, with some $0<\pi<1$. Then, letting
$\Theta_i=\{\theta_i\}$, $i=1,2$, in (\ref{204}), we have that
$$
W_1(\psi,\delta)=\pi P_{\theta_1}(\mbox{
reject}\;H_1)=\pi\alpha(\psi,\delta)$$ and
$$
W_2(\psi,\delta)=(1-\pi)P_{\theta_2}(\mbox{accept}\;
H_1)=(1-\pi)\beta(\psi,\delta),
$$
where $\alpha(\psi,\delta)$ and $\beta(\psi,\delta)$ are the type I
and type II error probabilities. Thus, taking in (\ref{204})
$w_1=\pi\alpha$, $w_2=(1-\pi)\beta$, with some
$\alpha,\beta\in(0,1)$, we see that (\ref{204}) is equivalent to
\begin{equation}\label{205}
\alpha(\psi,\delta)\leq \alpha,\quad\mbox{and}\quad
\beta(\psi,\delta)\leq \beta.
\end{equation}
Let now $\pi_2(\{\theta_0\})=1$ and suppose that the observations
are i.i.d.. Then our
problem of minimizing $N(\psi)=N(\theta_0;\psi)$ under restrictions
 (\ref{205}) is the classical Wald and Wolfowitz problem of minimizing the
expected sample size (see \cite{WaldWolfowitz}). It is well known that its solution is given by
the sequential probability ratio test (SPRT), and that it minimizes
the expected sample size under the alternative hypothesis as well
(see \cite{WaldWolfowitz}, \cite{Lorden}).

On the other hand, if $\pi_2(\{\theta\})=1$ with
$\theta\not=\theta_0$ and $\theta\not=\theta_1$, we have the problem
known as the modified Kiefer-Weiss problem, the problem of
minimizing the expected sample size, under $\theta$, among all
sequential tests subject to (\ref{205}) (see \cite{Weiss},
\cite{KieferWeiss}). The
general structure of the optimal sequential test in this problem is
given by Lorden \cite{Lorden} for i.i.d. observations.

So, we see  that considering natural particular cases  of sequential
procedures subject to (\ref{204}) and using different choices of $\pi_1$ in (\ref{203}) and
$\pi_2$ in (\ref{201})  we extend known
problems for i.i.d. observations to the case of general
discrete-time stochastic processes.

The method we use in this article was originally developed for testing of
two hypotheses \cite{NovikovSQA}, then extended for multiple hypothesis
 testing problems \cite{NovikovKybernetika2009}.
 An extension of the same method for hypothesis testing problems when
 control variables are present can be found in \cite{NovikovKybernetika2009Control}.

 A more general, than used in this article, setting for Bayes-type decision problems,
 where both the cost of observations and the loss functions
 depend on the true value of the parameter and on the observations, is considered in \cite{NovikovIWSM2009}.

From this time on, our aim will be minimizing $N(\psi)$, defined by
(\ref{201}), in the class of sequential statistical procedures
subject to (\ref{204}).

In Section 2, we reduce the problem to an optimal stopping problem.
In Section 3, we give a solution to the optimal stopping problems
in the class of truncated stopping rules, and in Section 4 in some natural class of non-truncated stopping rules.
In particular, in
Section 4 we give a solution to the problem of minimizing $N(\psi)$ in the class of all statistical procedures
satisfying $W_i(\psi,\delta)\leq w_i$, $i=1,\dots,k$ (see Remark \ref{r8}).

\section{\normalsize REDUCTION TO AN OPTIMAL STOPPING PROBLEM}\label{s2}

In this section, the problem of  minimizing the average sample
number (\ref{201}) over all sequential procedures subject to
(\ref{204}) will be reduced to an optimal stopping problem. This is
a usual treatment of conditional problems  in sequential hypothesis
testing (see, for example, \cite{Berk}, \cite{Lorden},
\cite{Castillo}, \cite{Muller}). We will use the same ideas to treat
the general statistical decision problem described above.

Let us define the following Lagrange-multiplier
function:
\begin{equation}{\label{4}}
L(\psi,\delta)=L(\psi,\delta;\lambda_1,\dots,\lambda_k)=N(\psi)+\sum_{i=1}^k\lambda_i W_i(\psi,\delta)
\end{equation}
where  $\lambda_i\geq 0$, $i=1,\dots, k$  are some constant
multipliers.

Let ${\Delta}$ be a class of sequential statistical procedures.

The following Theorem is a direct application of the method of
Lagrange multipliers  to the above optimization problem.

\begin{theorem}\label{t1} Let there exist $\lambda_i> 0$, $i=1,\dots,k$,
and a procedure $(\psi^*,\delta^*)\in {\Delta}$ such that for any
procedure $(\psi,\delta)\in {\Delta}$
\begin{equation}\label{5}
L(\psi^*,\delta^*;\lambda_1,\dots,\lambda_k)\leq L(\psi,\delta;\lambda_1,\dots,\lambda_k)
 \end{equation}
 holds and such that
 \begin{equation}\label{6}W_i(\psi^*,\delta^*)=w_i,\quad i=1,\dots
 k.
 \end{equation}

 Then for any test $(\psi,\delta)\in{\Delta}$ satisfying
 \begin{equation}\label{5bis}
 W_i(\psi,\delta)\leq w_i, \quad i=1,2,\dots ,k,
 \end{equation}
 it holds
\begin{equation}\label{5a}
N(\psi^*)\leq  N(\psi). \end{equation}

 The inequality
in (\ref{5a}) is strict if at least one of the inequalities
(\ref{5bis}) is strict.
\end{theorem}

{\bf Proof.}
 Let $(\psi,\delta)\in {\Delta}$ be any procedure
satisfying (\ref{5bis}). Because of (\ref{5}),
\begin{eqnarray}\label{5b}
L(\psi^*,\delta^*;\lambda_1,\dots,\lambda_k)=N(\psi^*)+\sum_{i=1}^k\lambda_i
W_i(\psi^*,\delta^*)\leq L(\psi,\delta;\lambda_1,\dots,\lambda_k)&&
\end{eqnarray}
\begin{equation}\label{5c}=N(\psi)+\sum_{i=1}^k\lambda_i W_i(\psi,\delta)\leq
N(\psi)+\sum_{i=1}^k\lambda_i w_i,
\end{equation}
%
where to get the last inequality we used (\ref{5bis}).
Taking into account conditions (\ref{6}) we get from this that
$$N(\psi^*)\leq N(\psi).
$$

To get the last statement of the theorem we note that if
$N(\psi^*)=N(\psi)$ then there are equalities in (\ref{5b}) --
(\ref{5c}) instead of the inequalities, which is only possible if
$W_i(\psi,\phi)=w_i$ for any $i=1,\dots,k$. $\Box$
\begin{remark}\label{r0}\rm
It is easy to see that defining a new loss function
$w^\prime(\theta,d)$ which is equal to $\lambda_i w(\theta,d)$
whenever $\theta\in\Theta_i$,$\quad i=1,\dots,k$, we have that the weighted
average loss $W(\psi,\delta)$ defined by (\ref{203}) with
$w(\theta,d)=w^\prime(\theta,d)$ coincides with the second summand
in (\ref{4}).

Because of this, we treat in what follows only the case of one
summand ($k=1$) in (\ref{4}),  being the Lagrange-multiplier
function defined as
\begin{equation}\label{208}
L(\psi,\delta;\lambda)=N(\psi)+\lambda W(\psi,\delta).
\end{equation}

It is obvious that the problem of minimization of (\ref{208}) is equivalent to that of minimization of
\begin{equation}\label{211}
R(\psi,\delta;c)=cN(\psi)+W(\psi,\delta),
\end{equation}
where $c>0$ is any constant, and, in the rest of the article, we will solve the problem of minimizing (\ref{211}), instead of (\ref{208}). This is because
the problem of minimization of (\ref{211}) is interesting by
itself, without its relation to the conditional problem above. For example, if $\pi_2=\pi_1=\pi$, it is easy to
see that it is equivalent to the problem of Bayesian sequential
decision-making, with the prior distribution $\pi$ and a fixed
cost $c$ per observation. The latter set-up is fundamental in
the sequential analysis (see \cite{Wald}, \cite{Ferguson},
\cite{DeGroot}, \cite{Zacks}, \cite{Ghosh}, 
among many others).
\end{remark}

 Because of Theorem \ref{t1},
from this time on, our main focus will be on the unrestricted
minimization of $R(\psi,\delta;c)$, over  all sequential decision procedures.

 Let us suppose, additionally to the assumptions of Introduction, that for any $n=1, 2\dots$ there exists a decision function
 $\delta_n^B=\delta_n^B(x_1,\dots,x_n)$ such that for any $d\in\mathscr D$
 \begin{equation}\label{206}
 \int w(\theta,d)f_\theta^n(x_1,\dots,x_n)d\pi_1(\theta)\geq \int
 w(\theta,\delta_n^B(x_1,\dots,x_n))f_\theta^n(x_1,\dots,x_n)d\pi_1(\theta)
 \end{equation}
 for $\mu^n$-almost all $(x_1,\dots,x_n)$.
Then $\delta_n^B$ is called the Bayesian decision function based on
$n$ observations. We do not discuss in this article the questions of
the existence of Bayesian decision functions, we just suppose that
they exist for any $n=1,2,\dots$ referring, e.g., to \cite{Wald} for an extensive underlying theory.

Let us denote by $l_n=l_n(x_1,\dots,x_n)$ the right-hand side of
(\ref{206}).
It easily follows from  (\ref{206}) that
\begin{equation}\label{2}
\int l_n d\mu^n=\inf_{\delta_n}\int E_\theta
w(\theta,\delta_n)d\pi_1(\theta),
\end{equation}
thus
$$
\int l_1 d\mu^1\geq \int l_2 d\mu^2\geq \dots.
$$

Because of that,  we suppose that
$$
\int l_1(x) d\mu(x)<\infty
$$
which makes all the Bayesian risks (\ref{2}) finite, for any
$n=1,2,\dots$.

Let $\delta^B=(\delta_1^B,\delta_2^B, \dots)$. The following Theorem
shows that the only decision rules worth our attention are the Bayesian ones. Its ``if''-part is, in essence,
 Theorem 5.2.1 \cite{Ghosh}.

Let for any $n=1,2,\dots$ and for any stopping rule $\psi$ $$s_n^\psi=(1-\psi_1)\dots(1-\psi_{n-1})\psi_n,
$$ and let $$ S_n^\psi=
\{(x_1,\dots,x_n): s_n^\psi(x_1,\dots,x_n)>0
\}
$$
for all $n=1,2,\dots$.

\begin{theorem} \label{t2} For
any sequential procedure $(\psi,\delta)$
\begin{equation}\label{6a}
W(\psi,\delta)\geq W(\psi,\delta^B)=\sum_{n=1}^\infty
\int s_n^\psi l_n d\mu^n.
\end{equation}

Supposing that the right-hand side of (\ref{6a}) is finite, the equality in (\ref{6a}) is only possible if
$$
\int
 w(\theta,\delta_n)f_\theta^nd\pi_1(\theta)=\int
 w(\theta,\delta_n^B)f_\theta^nd\pi_1(\theta)
$$
$\mu^n$-almost everywhere on $S_n^\psi$ for all $n=1,2,\dots$.
\end{theorem}
{\bf Proof.} It is easy to see that $W(\psi,\delta)$ on the left-hand side of (\ref{6a}) has the following equivalent form:
\begin{equation}\label{7f}
W(\psi,\delta)=\sum_{n=1}^\infty \int s_n^\psi
\int w(\theta,\delta_n)f_\theta^nd\pi_1(\theta) d\mu^n.
 \end{equation}

Applying (\ref{206}) under the integral sign in each summand in (\ref{7f}) we immediately
have:
 \begin{equation}\label{7g}
W(\psi,\delta)\geq \sum_{n=1}^\infty \int s_n^\psi
\int w(\theta,\delta_n^B)f_\theta^nd\pi_1(\theta) d\mu^n=W(\psi,\delta^B).
 \end{equation}
If $W(\psi,\delta^B)<\infty$, then (\ref{7g}) is equivalent to
$$
 \sum_{n=1}^\infty \int s_n^\psi\Delta_n
 d\mu^n\geq 0,
$$
where
$$\Delta_n=\int w(\theta,\delta_n)f_\theta^nd\pi_1(\theta)-\int w(\theta,\delta_n^B)f_\theta^nd\pi_1(\theta),$$
which is, due to (\ref{206}), non-negative $\mu^n$-almost everywhere for all $n=1,2,\dots$.
Thus, there is an equality in (\ref{7g}) if and only if $\Delta_n=0$ $\mu^n$-almost everywhere on $S_n^\psi=\{s_n^\psi>0\}$ for all $n=1,2,\dots$.
$\Box$

Because of (\ref{211}), it follows from Theorem \ref{t2} that for any sequential
decision procedure $(\psi,\delta)$
\begin{equation}\label{14}
R(\psi,\delta;c)\geq R(\psi,\delta^B;c).
\end{equation}
The following lemma gives the right-hand side of (\ref{14}) a more
convenient form.

For any probability measure $\pi$ on $\Theta$ let us denote
$$
P^{\pi}(\tau_\psi=n)\equiv\int
P_\theta(\tau_\psi=n)d\pi(\theta)=\int E_\theta
s_n^\psi d\pi(\theta),
$$
for $n=1,2,\dots$ Respectively,
$P^{\pi}(\tau_\psi<\infty)=\sum_{n=1}^\infty
P^{\pi}(\tau_\psi=n)$, and $$E^{\pi}\tau_\psi=\int E_\theta\tau_\psi d\pi(\theta).$$
\begin{lemma}\label{l0}
If
\begin{equation}\label{3}
P^{\pi_2}(\tau_\psi<\infty)=1
\end{equation}
then
\begin{equation}\label{7}
R(\psi,\delta^B;c)=\sum_{n=1}^\infty \int s_n^\psi \left(cnf^n+l_n \right)d\mu^n,
\end{equation}
where, by definition,
\begin{equation}\label{207}
f^n=f^n(x_1,\dots,x_n)=\int f_\theta^n(x_1,\dots,x_n)d\pi_2(\theta).
\end{equation}
\end{lemma}
{\bf Proof.} By Theorem \ref{t2},
\begin{equation}\label{12}
R(\psi,\delta^B;c)=cN(\psi)+W(\psi,\delta^B)=
cN(\psi)+\sum_{n=1}^\infty \int s_n^\psi
l_n d\mu^n.
\end{equation}
If now (\ref{3}) is fulfilled, then, by the Fubini theorem,
$$
N(\psi)=\int \sum_{n=1}^\infty n
E_\theta s_n^\psi d\pi_2(\theta)=\sum_{n=1}^\infty \int
E_\theta n s_n^\psi d\pi_2(\theta)
$$
$$
=\sum_{n=1}^\infty \int
s_n^\psi\left(n\int f_\theta^n
d\pi_2(\theta)\right)d\mu^n=\sum_{n=1}^\infty \int
s_n^\psi nf^n d\mu^n,
$$
so, combining this with (\ref{12}), we get (\ref{7}). $\Box$

Let us denote
\begin{equation}\label{15}
R(\psi)=R(\psi;c)=R(\psi,\delta^B;c).
\end{equation}
By Lemma \ref{l0},
\begin{equation}\label{16}
R(\psi)=\begin{cases}\displaystyle\sum_{n=1}^\infty \int
s_n^\psi \left(cnf^n+l_n
\right)d\mu^n,\;\mbox{if}\; P^{\pi_2}(\tau_\psi<\infty)=1,\cr
\infty,\;\mbox{otherwise}.
\end{cases}
\end{equation}

The aim of what follows is to minimize $R(\psi)$ over all stopping
rules. In this way, our problem of minimization of $R(\psi,\delta)$
is reduced to an optimal stopping problem.

\section{\normalsize OPTIMAL TRUNCATED STOPPING RULES}
In this section, as a first step, we characterize the structure of optimal stopping rules in the class ${\mathscr F}^N$, $N\geq 2$, of all
truncated stopping rules, i.e., such that
\begin{equation}\label{11}
\psi=(\psi_1,\psi_2,\dots,\psi_{N-1},1,\dots)\end{equation}
(if $(1-\psi_1)\dots (1-\psi_n)=0$ $\mu^n$-almost everywhere for some $n<N$, we suppose that $\psi_k\equiv 1$ for any $k>n$, so ${\mathscr F}^N\subset {\mathscr F}^{N+1}$, $N=1,2,\dots$).

Obviously, for any $\psi\in{\mathscr F}^N$
$$R(\psi)=R_N(\psi)=\sum_{n=1}^{N-1}\int
s_n^\psi(cnf^n+l_n)d\mu^n
+\int
t_N^\psi\left(cNf^N+l_N\right)d\mu^N,
$$
where for any $n=1,2,\dots$ $$t_n^\psi=t_n^\psi(x_1,\dots,x_n)=(1-\psi_1(x_1))(1-\psi(x_1,x_2))\dots
(1-\psi_{n-1}(x_1,\dots,x_{n-1}))$$
(we suppose, by definition, that $t_1^\psi\equiv 1$).

Let us introduce a sequence of functions $V_n^N$, $n=1,\dots, N$, which will define optimal stoppings rules. Let $V_N^N\equiv l_N$, and recursively for $n=N-1, N-2, \dots 1 $
\begin{equation}\label{48}
V_n^N=\min\{l_n,Q_n^N\},
\end{equation}
where \begin{equation}\label{48a}Q_n^N=Q_n^N(x_1,\dots,x_n)=cf^n(x_1,\dots,x_n)+\int
V_{n+1}^N(x_1,\dots,x_{n+1})d\mu(x_{n+1}),\end{equation}
$n=0,1,\dots,N-1$ (we assume that $f^0\equiv 1$). Please, remember that all $V_n^N$ and $Q_n^N$ implicitly depend on the ``unitary observation cost'' $c$.

The following theorem characterizes the structure of optimal stopping rules in ${\mathscr F}^N$.
\begin{theorem}\label{t3} For all $\psi\in{\mathscr F}^N$
\begin{equation}\label{46b}
    R_N(\psi)\geq Q_0^N.
\end{equation}

 The lower bound in (\ref{46b}) is attained by a $\psi\in{\mathscr F}^N$ if and only if
\begin{eqnarray}\label{49}
I_{\{l_{n}< Q_n^N \}}\leq \psi_{n}\leq I_{\{l_{n}\leq
Q_n^N \}}
\end{eqnarray}
$\mu^n$-almost everywhere on $$T_n^\psi=\{(x_1,\dots,x_n):t_n^\psi(x_1,\dots,x_n)>0\},$$ for all $n=1,2,\dots,
N-1$.

\end{theorem}

The proof of Theorem \ref{t3} can be conducted following the lines of the proof of Theorem 3.1 in
 \cite{NovikovSQA} (in a less formal way, the same routine is used to obtain
 Theorem 4 in \cite{NovikovKybernetika2009}).
 In fact, both of these theorems are particular cases of Theorem \ref{t3}.
\begin{remark}\label{r3a}\rm Despite that $\psi$ satisfying (\ref{49}) is optimal among
all truncated  stopping rules in $\mathscr F^N$, it only makes practical sense if
\begin{equation}\label{31}l_0=\inf_{d}\int w(\theta,d)d\pi_1(\theta)\geq Q_0^N.\end{equation}

 Indeed, if (\ref{31}) does not hold, we can, without taking any observation, make any decision $d_0$ such that
 $\int w(\theta,d_0)d\pi_1(\theta)<Q_0^N$, and this guarantees that
 this trivial procedure (something like ``$(\psi_0,d_0)$'' with $R(\psi_0,d_0)=\int w(\theta,d_0)d\pi_1(\theta)<Q_0^N$)
 performs better than the best procedure with the optimal stopping time in $\mathscr F^N$.

Because of this, $V_0^N$, defined by (\ref{48}) for
$n=0$, may be considered the  ``minimum
value of $R(\psi)$'', when taking no observations is allowed.\\
\end{remark}
\begin{remark}\label{r4} \rm
When $\pi_2$ in (\ref{201}) coincides with $\pi_1$ in (\ref{203}) (Bayesian setting),
an optimal truncated (non-randomized) stopping rule for minimizing
(\ref{211}) is provided by Theorem 5.2.2 in \cite{Ghosh}. Theorem \ref{t3} describes the class of {\em all randomized} optimal
stopping rules for the same problem in this particular case. This
may be irrelevant if one is interested in the purely Bayesian
problem, because any of these stopping rules provides the same
minimum value  of the risk.

Nevertheless, this extension of the class of optimal procedures may
be useful for complying with (\ref{6}) in Theorem \ref{t1} when
seeking for optimal sequential procedures for the original
conditional problem (minimization of $N(\psi)$ given that
$W_i(\psi,\delta)\leq w_i$, $i=1,\dots, k$, see Introduction and the
discussion therein). This is very much like in non-sequential
hypothesis testing, where the randomization is crucial for finding
the optimal level-$\alpha$ test in the Neyman-Pearson problem (see,
for example, \cite{Lehmann}).
\end{remark}

\section{\normalsize OPTIMAL NON-TRUNCATED STOPPING RULES}
In this section, we solve the problem of minimization of $R(\psi)$ in natural classes of non-truncated stopping rules $\psi$.

Let $\psi$ be any stopping rule. Define
\begin{equation}
\label{50a}
R_N(\psi)=R_N(\psi;c)=\sum_{n=1}^{N-1}\int
s_n^\psi(cnf^n+l_n)d\mu^n+\int t_N^\psi\left(cNf^N+l_N\right)d\mu^{N}.
\end{equation}
This is the ``risk'' (\ref{211})  for $\psi$ truncated  at
$N$, i.e. the rule with the components
$\psi^N=(\psi_1,\psi_2,\dots,\psi_{N-1},1,\dots)$:
$R_N(\psi)=R(\psi^N)$.

Because $\psi^N$ is truncated, the results of the preceding section
apply, in particular, the lower bound of (\ref{46b}). Very much like in \cite{NovikovSQA} and in \cite{NovikovKybernetika2009}, our aim is to pass to the limit, as $N\to\infty$, in order to obtain a lower bound for $R(\psi)$, and conditions for attaining this bound.

It is easy to see that $V_n^N(x_1,\dots,x_n)\geq V_n^{N+1}(x_1,\dots,x_n)$ for all $N\geq n$, and for all $(x_1,\dots,x_n)$, $n\geq 1$ (see, for example, Lemma  3.3 in \cite{NovikovSQA}). Thus, for any $n\geq 1$ there exists
$$
V_n=V_n(x_1,\dots,x_n)=\lim_{N\to\infty} V_n^N(x_1,\dots,x_n),
$$
($V_n$ implicitly depend on $c$, as $V_n^N$ do). It immediately follows from the monotone convergence theorem that for all $n\geq 1$
\begin{equation}\label{50b}
\lim_{N\to\infty}Q_n^N(x_1,\dots,x_n)=cf^n(x_1,\dots,x_n)+\int
V_{n+1}(x_1,\dots,x_{n+1})d\mu(x_{n+1})
\end{equation}
(see (\ref{48a})).
Let $Q_n=Q_n(x_1,\dots,x_n)=\lim_{N\to\infty}Q_n^N(x_1,\dots,x_n)$.

In addition, passing to the limit, as $N\to\infty$, in (\ref{48}) we obtain
$$
V_n=\min\{l_n,Q_n\},\quad n=1,2,\dots.
$$

Let now $\mathscr F$ be any class of stopping rules such that $\psi\in\mathscr F$ entails $R_N(\psi)\to R(\psi)$,
as $N\to\infty$.  It is easy to see that such classes exist, for example, any $\mathscr F^N$ has this property.
Moreover, we will assume that all truncated stopping rules are included in $\mathscr F$,
i.e. that $\bigcup_{N\geq 1}\mathscr F^N\subset \mathscr F$.

It follows from Theorem \ref{t3} now that for all $\psi\in \mathscr F$
\begin{equation}\label{50c}
R(\psi)\geq Q_0.
\end{equation}
The following lemma states that, in fact, the lower bound in (\ref{50c}) is the infimum of the risk $R(\psi)$ over $\psi\in\mathscr F$.
\begin{lemma}\label{l1}
$$Q_0=\inf_{\psi\in\mathscr F}R(\psi).$$
\end{lemma}
The proof of Lemma \ref{l1} is very close to that of Lemma 3.5 in \cite{NovikovSQA} (see also Lemma 6 in \cite{NovikovKybernetika2009}) and is omitted here.
\begin{remark}\label{r5}\rm
Again (see Remark \ref{r4}), if $\pi_1=\pi_2$,  Lemma \ref{l1} is essentially
Theorem 5.2.3 in \cite{Ghosh} (see also  Section 7.2 of \cite{Ferguson}) .
 \end{remark}
The following Theorem gives the structure of optimal stopping rules in $\mathscr F$.

\begin{theorem}\label{t4}
If there exists $\psi\in\mathscr F$ such that
\begin{equation}\label{48b}
R(\psi)=\inf_{\psi^\prime\in\mathscr F}R(\psi^\prime),
\end{equation}
then
\begin{eqnarray}\label{49a}
I_{\{l_{n}< Q_n \}}\leq \psi_{n}\leq I_{\{l_{n}\leq
Q_n \}}
\end{eqnarray}
$\mu^n$-almost everywhere on $T_n^\psi$ for all $n=1,2,\dots$.

On the other hand, if a stopping rule $\psi$ satisfies (\ref{49a}) $\mu^n$-almost everywhere on $T_n^\psi$ for all $n=1,2,\dots$,
and $\psi\in \mathscr F$, then $\psi$ satisfies (\ref{48b}) as well.
\end{theorem}

The proof of Theorem \ref{t4} is very close to the proof of Theorem 3.2 in \cite{NovikovSQA} or Theorem 6 in
\cite{NovikovKybernetika2009} and is omitted here.

It follows from Theorem \ref{t4} that ``$\psi\in\mathscr F$''
is a sufficient condition for the optimality of a stopping rule $\psi$ satisfying (\ref{49a}).
In the hypothesis testing problems considered in \cite{NovikovSQA} and in \cite{NovikovKybernetika2009},
there are large classes of problems (called truncatable) for which $R_N(\psi)\to R(\psi)$, as $N\to\infty$, for {\em all}
stopping times $\psi$. In this article, we also identify the problems where this is the case.

Let us say that a stopping rule is {\em truncatable} if  $R_N(\psi)\to R(\psi)$, as $N\to\infty$.
It is obvious that all truncated stopping rules are truncatable. In particular, Theorem \ref{t4} holds when
$\mathscr F$ is the set of all truncatable stopping rules.

The following Lemma gives a necessary and sufficient condition for truncatablity of a stopping rule.
\begin{lemma}\label{l2} A stopping rule $\psi$ with $R(\psi)<\infty$ is truncatable
if and only if
\begin{equation}\label{13}
\int t_N^\psi l_Nd\mu^N\to
0,\quad\mbox{as}\quad N\to\infty.
\end{equation}
If $R(\psi)=\infty$, then $R_N(\psi)\to \infty$, $N\to\infty$.
\end{lemma}
{\bf Proof.}
Let $\psi$ be such that $R(\psi)<\infty$. 

Suppose that (\ref{13}) is fulfilled.
 Then, by (\ref{50a})
 \begin{equation}\label{49d}
 R(\psi)-R_N(\psi)=
\sum_{n=N}^\infty \int
c_n^\psi(cnf^n+l_n)d\mu^n -c\int
t_N^\psi Nf^Nd\mu^N+\int t_N^\psi l_N d\mu^N.
\end{equation}

The first summand converges to zero, as $N\to\infty$, being the tail
of a convergent series (this is because $R(\psi)<\infty$).

The third summand in (\ref{49d}) goes to 0 as $N\to\infty$, because of (\ref{13}).

The integral in the second summand in (\ref{49d}) is equal to
\begin{equation*}
N
P^{\pi_2}(\tau_\psi\geq N)\leq E^{\pi_2}\tau_\psi I_{\{\tau_\psi\geq N\}}\to 0,
\end{equation*}
as $N\to\infty$, because $E^{\pi_2}\tau_\psi<\infty$ (this is due to $R(\psi)<\infty$ again).

It follows from (\ref{49d}) now that $R_N(\psi)\to R(\psi)<\infty$ as $N\to\infty$.

Let us suppose now that $R_N(\psi)\to R(\psi)<\infty$ as $N\to\infty$. For the same reasons as above,
the first two summands on the right-hand side of
(\ref{49d}) tend to 0 as $N\to\infty$, therefore so does the third, i.e. (\ref{13}) follows.
The first assertion of Lemma \ref{l2} is proved.

If  $R(\psi)=\infty$,
this may be because $P^{\pi_2}(\tau_\psi<\infty)<1$, or, if not,
because
$$
\sum_{n=1}^{\infty}\int s_n^\psi(cnf^n+l_n)d\mu^n=\infty.
$$
In the latter case, obviously, $$R_N(\psi)\geq \sum_{n=1}^{N-1}\int
s_n^\psi(cnf^n+l_n)d\mu^n\to\infty,\quad
\mbox{as}\quad N\to\infty.$$

In the former case,

$$R_N(\psi)\geq c\int t_N^\psi Nf^N d\mu^N=cNP^{\pi_2}(\tau_\psi\geq N)\to\infty,$$
as $N\to\infty$, as well.
$\Box$

Let us say that the problem (of minimization of $R(\psi)$) is truncatable if all stopping rule $\psi$ are
truncatable.

Corollary  \ref{c1} below gives some practical sufficient conditions for truncatability of a problem.

\begin{corollary}\label{c1}
The problem of minimization of $R(\psi)$ is truncatable if\\
 {\rm i)} the loss function $w$ is
 bounded, and
 \begin{equation}\label{19}R(\psi)<\infty\quad\mbox{implies that}\quad P^{\pi_1}(\tau_\psi<\infty)=1,\end{equation}
 or\\
{\rm ii)} \begin{equation}\label{18}\int l_N d\mu^N\to 0,\end{equation}
as $N\to\infty$.
\end{corollary}
{\bf Proof.} If $w(\theta,d)<M<\infty$ for any $\theta$ and $d$,  then, by the
 definition of $l_N$,
 \begin{equation}\label{20}
 \int t_N^\psi l_N d\mu^N\leq M \int t_N^\psi\left( \int
 f_\theta^N d\pi_1(\theta)\right) d\mu^N=MP^{\pi_1}(\tau_\psi\geq N).
 \end{equation}
 If now $R(\psi)<\infty$, then by (\ref{19})  the right-hand side of (\ref{20}) tends to 0,  as $N\to\infty$,
 i.e. (\ref{13}) is fulfilled for any $\psi$ such that $R(\psi)<\infty$. Thus, by Lemma \ref{l2} any $\psi$ is truncatable.

If (\ref{18}) is fulfilled, then (\ref{13}) is satisfied for any
$\psi$. Again, by Lemma \ref{l2} any $\psi$ is truncatable. $\Box$

\begin{remark}\label{r1} \rm
Condition {\rm i)} of Corollary \ref{c1} is fulfilled for any Bayesian hypothesis testing
problem (i.e. when $\pi_1=\pi_2=\pi$) with bounded loss function (see, for example,
 \cite{NovikovSQA} and \cite{NovikovKybernetika2009}). Indeed, in this case $R(\psi)<\infty$ implies
 $E^\pi \tau_\psi<\infty$, so, in particular,  $P^\pi(\tau_\psi<\infty)=1$.
\end{remark}
\begin{remark}\label{r2}\rm
It is easy to see that Condition {\rm ii)} of Corollary \ref{c1} is equivalent to
$$
\int_\Theta E_\theta w(\theta,\delta_{N}^B)d\pi_1(\theta)\to 0, \;N\to\infty,
$$
i.e. that the Bayesian risk, with respect to the prior distribution $\pi_1$, of an optimal
procedure based on sample of a fixed size $N$, vanishes as $N\to\infty$. This is a very typical
behavior of statistical risks.
\end{remark}

The following Theorem is an immediate consequence of Theorem \ref{t4}.
\begin{theorem}\label{t5} Let the problem of minimization of $R(\psi)$ be truncatable, and let $\mathscr F$ be the set of
all stopping rules.
Then
\begin{equation}\label{48c}
R(\psi)=\inf_{\psi^\prime\in\mathscr F}R(\psi^\prime)
\end{equation}
if and only if
\begin{eqnarray}\label{49b}
I_{\{l_{n}< Q_n \}}\leq \psi_{n}\leq I_{\{l_{n}\leq
Q_n \}}
\end{eqnarray}
$\mu^n$-almost everywhere on $T_n^\psi$ for all $n=1,2,\dots$.
\end{theorem}
\begin{remark}\label{r7}\rm
Once again (see Remark \ref{r3a}), the optimal stopping rule $\psi$
from Theorem \ref{t5} (and Theorem \ref{t4}) only makes practical sense if $l_0\geq Q_0=\inf_{\psi\in\mathscr F}R(\psi)$,
because otherwise the trivial rule, which does not take
any observation, performs better than $\psi$, from the point of view of minimization of $R(\psi)$.
\end{remark}
\begin{remark}\label{r8}\rm
Combining Theorems \ref{t1}, \ref{t2}, and \ref{t5} we
immediately have the following solution to the conditional problem posed in Introduction.

Let $\lambda_1,\dots$, $\lambda_k$ be arbitrary positive constants.
Let $\delta_n^B$, $n=1,2,\dots$ be Bayesian, with respect to $\pi_1$, decision rules for the ``loss function''
$$
w^\prime(\theta,d)=\sum_{i=1}^k\lambda_iw(\theta,d)I_{\Theta_i}(\theta),
$$
i.e. such that for all $d\in\mathscr D$
\begin{equation}\label{60}
\sum_{i=1}^k\lambda_i\int_{\Theta_i} w(\theta,d)f_\theta^nd\pi_1(\theta)
\geq l_n=\sum_{i=1}^k\lambda_i\int_{\Theta_i} w(\theta,\delta_n^B)f_\theta^nd\pi_1(\theta)
\end{equation}
$\mu^n$-almost everywhere (remember that $\delta_n^B=\delta_n^B(x_1,\dots,x_n)$
and $f_\theta^n=f_\theta^n(x_1,\dots,x_n)$).

For any $N\geq 1$ define
$
V_N^N=l_N,$ and $ V_n^N=\min\{l_n,Q_n^N\}$ for $ n=N-1,N-2,\dots,1$,
where $Q_n^N=f^n+\int l_{n+1}d\mu(x_{n+1})$, with $f^n=\int_\Theta f_\theta^nd\pi_2(\theta)$.

Let also $V_n=\lim_{N\to\infty}V_n^N$ and $Q_n=\lim_{N\to\infty}Q_n^N$, $n=1,2,\dots$.

Suppose, finally, that
the problem is truncatable (see Corollary \ref{c1} for sufficient conditions for that).

Let $\psi$ be any stopping rule satisfying
\begin{equation}
I_{\{l_{n}< Q_n \}}\leq \psi_{n}\leq I_{\{l_{n}\leq
Q_n \}}
\end{equation}
$\mu^n$-almost everywhere on $T_n^\psi$ for all $n=1,2,\dots$.

 Then for any sequential decision procedure
$(\psi^\prime,\delta)$ such that
\begin{equation}\label{25}
W_i(\psi^\prime,\delta)\leq W_i(\psi,\delta^B),\quad i=1,\dots, k,
\end{equation}
it holds
\begin{equation}\label{26}
N(\psi)\leq N(\psi^\prime).
\end{equation}
The inequality in (\ref{26}) is strict if at least one of the
inequalities in (\ref{25}) is strict.

If there are equalities in all of the inequalities in (\ref{25}) and
(\ref{26}), then \begin{equation}
I_{\{l_{n}< Q_n \}}\leq \psi_{n}^\prime\leq I_{\{l_{n}\leq
Q_n \}}
\end{equation}
$\mu^n$-almost everywhere on $T_n^{\psi^\prime}$ for all $n=1,2,\dots$, and
$$
\sum_{i=1}^k\lambda_i\int_{\Theta_i} w(\theta,\delta_n)f_\theta^nd\pi_1(\theta)=\sum_{i=1}^k\lambda_i\int_{\Theta_i} w(\theta,\delta_n^B)f_\theta^nd\pi_1(\theta)
$$
$\mu^n$-almost everywhere on $S_n^{\psi^\prime}$ for all $n=1,2,\dots$.
\end{remark}
For Bayesian problems (when $\pi_1=\pi_2=\pi$) Theorem \ref{t5} can be reformulated in the following equivalent way.

Let $$R_n=\frac{l_n}{f^n}=\frac{\int_\Theta f_\theta^n w(\theta,\delta_n^B)d\pi(\theta)}
{\int_\Theta f_\theta^n d\pi(\theta)}$$
be the posterior risk (see, e.g., \cite{Berger}). Let $v_N^N\equiv R_N(X_1,\dots,X_n)$, and recursively for $n=N-1,N-2,\dots,1$
$$
v_n^N(X_1,\dots,X_n)=\min\{R_n(X_1,\dots,X_n),q_n^N(X_1,\dots,X_n) \},
$$
where $$q_n^N(X_1,\dots,X_n)=c+E^\pi \{v_{n+1}^N|X_1,\dots,X_n\}$$ ($E^\pi$ stands for the expectation with respect
to the  family of finite-dimensional densities $f^n=\int_\Theta f_\theta^n d\pi(\theta)$, $n=1,2,\dots$, meaning, in particular, that
$$E^\pi \{v_{n+1}^N|x_1,\dots,x_n\}=\int \frac{v_{n+1}^N(x_1,\dots,x_{n+1})f^{n+1}(x_1,\dots,x_{n+1})}{f^{n}(x_1,\dots,x_{n})}d\mu(x_{n+1})).$$

Let, finally, $v_n=v_n(X_1,\dots,X_n)=\lim_{N\to\infty}v_n^N(X_1,\dots,X_n)$,
and \\$q_n=q_n(X_1,\dots,X_n)$ $=\lim_{N\to\infty}q_n^N(X_1,\dots,X_n)$, $n=1,2,\dots$.

Then, the following reformulation of Theorem \ref{t5} gives, for a truncatable Bayesian problem, the
structure of all  Bayesian randomized tests (cf. Theorem 7, Ch. 7, in \cite{Berger}).

\begin{theorem}\label{t6} Let the problem of minimization of $R(\psi)$ be truncatable, and let $\mathscr F$ be the set of
all stopping rules.
Then
\begin{equation}\label{48d}
R(\psi)=\inf_{\psi^\prime\in\mathscr F}R(\psi^\prime)
\end{equation}
if and only if
\begin{eqnarray}\label{49c}
I_{\{R_{n}< q_n \}}\leq \psi_{n}\leq I_{\{R_{n}\leq
q_n \}}
\end{eqnarray}
$P^\pi$-almost surely on $T_n^\psi$ for all $n=1,2,\dots$.
\end{theorem}

\begin{remark}\label{r3}\rm
More general variants of Theorem \ref{t6}, for cases when the loss function due
to incorrect decision is of the form $w(\theta,d)=w_n(\theta,d;x_1,\dots,x_n)$ and/or the cost
of the observations $(x_1,\dots,x_n)$
is of type $K_\theta^{n}(x_1,\dots,x_{n})$, can  easily be  deduced from  Theorem 4 \cite{NovikovIWSM2009}.
In particular, this gives the structure of optimal sequential multiple hypotheses tests
for the problem considered in
Section 9.4 of \cite{Zacks}.
\end{remark}
\begin{remark}\label{r6}\rm
Theorem \ref{t6}, in particular, gives a solution to optimal sequential hypothesis testing problems considered in
\cite{CochlarVrana} and \cite{Cochlar} (where the general theory of optimal stopping is used, see \cite{Chow} or \cite{Shiryaev}).
See \cite{NovikovSQA} and \cite{NovikovKybernetika2009} for a more detailed description of  the respective Bayesian sequential procedures.
\end{remark}

\normalsize
\section*{\normalsize ACKNOWLEDGEMENTS}
\small
The author thanks the anonymous referees for their valuable comments
and suggestions on the improvement of an early version of this paper.

 The author greatly appreciates the support of the Autonomous
 Metropolitan University, Mexico City, Mexico, where this work was
 done, and the support of the National System of Investigators (SNI)
 of
 CONACyT, Mexico.

This work is also partially supported by Mexico's CONACyT Grant no.
CB-2005-C01-49854-F. \vspace{3mm}
 \footnotesize
\begin{flushright}
(Received \today)\,\ \rule{0mm}{0mm}
\end{flushright}

\vspace*{2mm}

{\mi
\begin{flushright}
\begin{minipage}[]{124mm}
{Andrey Novikov, Departamento de Matem\'aticas, Universidad
Aut\'onoma Metropolitana - Unidad Iztapalapa, San Rafael Atlixco
186, col. Vicentina, C.P. 09340, M\'exico D.F., M\'exico
\\ e-mail: {\tt an@xanum.uam.mx}
\\{\tt http://mat.izt.uam.mx/profs/anovikov/en}}
\end{minipage}
\end{flushright}
}

\end{document}